\documentclass{article}
\usepackage{amssymb,amsmath}
\usepackage{a4wide}
\usepackage{theorem,latexsym,natbib}

\usepackage{pdfsync}
\usepackage{hyperref,natbib}
\newcommand\BE{\ensuremath{\mathbb{E}}}

\newcommand\ind{\ensuremath{\mathbf{1}}}

\newcommand\BN{\ensuremath{\mathbb{N}}}

\def\Var{\mathop{\rm Var}\nolimits}
\def\Abs#1{\left\vert #1 \right\vert}  
\def\Re{\textrm{Re\,}}

\def\d{\textrm{\,d}}

\newtheorem{Theo}{Theorem}[section]
\newtheorem{Lem}{Lemma}[section]

\newtheorem{Prop}{Proposition}[section]
{\theorembodyfont{\rmfamily}\newtheorem{Rem}{Remark}}
{\theorembodyfont{\rmfamily}}
\newtheorem{Def}{Definition}[section]

\title{Marginal density estimation for linear processes with
cyclical long memory } 
\author {Mohamedou OULD HAYE${}^1$ and Anne PHILIPPE${}^2$} 
\date{ \small ${}^1$School of Mathematics and Statistics, Carleton
 University \\ 1125 Colonel By Drive \\ Ottawa, Ontario K1S 5B6,
 Canada\\ 
 ${}^2$ Laboratoire de Mathématiques Jean Leray
UMR CNRS 6629\\ Universit\'e de Nantes \\ 
2 rue de la Houssini\`ere - BP 92208
44322 Nantes Cedex 3, France
 }

\begin{document}
\maketitle

\begin{center} {\bf Abstract}
\end{center}

Some convergence results on the kernel density estimator are proven
for a class of linear processes with  cyclic effects. In particular
we extend the results of \cite{0862.60026} and  \cite{MR1457496,0695.60043}  to
the stationary processes for  which the singularities of the spectral density are not limited to the origin. We show that the 
convergence rates and the limiting distribution may be different in this context.


\setcounter{equation}{0}
\section{Introduction}

\cite{Hosking-81} introduced long memory processes with quasi 
periodic behaviour. This fact corresponds, for stationary processes,  to spectral densities which exhibit 
singularities at non zero frequencies. Many authors have contributed to the construction of 
fractional models with singularities/poles outside the origin, see 
for instance,
\cite{Gray:1994,Gray:1989,Hassler-94,viano95,leipus:viano:00, MR2532092}. 

We can distinguish between two types of long memory: one regular and
the other  cyclical  according to whether the spectral density has a pole at
the origin or outside the origin. From a statistical point of view, the estimators of the long memory parameter 
have been adapted to yield some estimates if cyclical effects are
assumed. In a parametric context, the $\sqrt{n}$-consistency of the maximum
likelihood estimate or the Whittle estimate has been proved (see
\cite{Hosoya-97,Giraitis:Hidalgo:Robinson-01} when the
pole is unknown). Semi parametric estimates can be more or less easily adapted
to the \textit{cyclical} case (see \cite{1051.62075,0974.62079,arteche:robinson,MR2497555,MR2201232,MR2060018}). 

When we consider empirical process related statistics, the situation is more delicate. The normalisation and the limit
distribution can be different according to  whether the memory is regular
or cyclical. An important literature is devoted to the convergence of
the empirical process, see for instance \cite{0862.60026,MR1713796} in regular case and \cite{MR1943152}
\cite{ouldhaye:aphil:2003} in cyclical case.

In this paper we give some convergence results on the kernel
estimator of the marginal density $f$. 
 Let $(X_1,\cdots,X_n)$ be an observed sample from $f$, 
the kernel estimator of $f$ is defined by 
\begin{equation}\label{noyau}
\tilde{f}_n(x)=\frac{1}{nm_n}\sum_{j=1}^nK\bigl(\frac{x-X_j}{m_n}\bigr). 
\end{equation}
where the bandwidth  $m_n$  is a sequence such that $m_n \to +\infty$
and $nm_n \to 0$ as $n\to\infty$, and $K$ is a 
kernel function.

Consider the following infinite moving average process,
\begin{equation}\label{lineaires}
X_t=\sum_{j=-\infty}^{t}b_{t-j}\xi_j , \qquad t\ge 1
\end{equation}
where 
\begin{itemize}
\item the sequence $(b_k)_k$ has the form 
\begin{equation}\label{b(s)}
b_k=k^{-(\alpha+1)/2}\sum_{j\in J}a_j\bigl(\cos
k\lambda_j+o(1)\bigr), \qquad k\to\infty
\end{equation}
where $\alpha \in (0,1)$ and $\lambda_j \not= 0 $ for all $j\in J$, a
finite non empty subset of $\BN$. 
\item $(\xi_n)_n$ is a sequence of 
independent and identically distributed random variables  with zero mean
and finite variance $\BE \xi_0^2=\sigma^2< \infty$. 
\end{itemize}

From \cite{Giraitis:leipus:1995}, the covariance
function $r$ of $(X_t)$ defined by  \eqref{lineaires} and \eqref{b(s)} has 
the form
\begin{equation}\label{r(n)}
r(h)=h^{-\alpha}\sum_{j\in J}a_j\bigl(\cos h\lambda_j+o(1)\bigr). 
\end{equation}
as $h$ tends to infinity.

A large class of linear processes satisfying these
conditions is obtained by filtering a white noise $(\xi_i)$  
\begin{equation}\label{g1}
X_t = G(B)\xi_t \quad \mathrm{ with }\quad G(z)=g(z)\prod_{j=-m}^m\bigl(1-e^{i\lambda_j}z
\bigr)^{(\alpha_j-1)/2},\qquad m\ge1,
\end{equation}
where $B$ is the backshift operator and 
where $g$ is an analytic function on $\{ |z|< 1 \}$,
 continuous on $\{ |z|\le 1 \}$ and
$g(z)\not= 0$ if $\Abs{z}=1$, 
and where  
$$
0<\alpha_j\leq 1,\quad\alpha_j=\alpha_{-j},\quad\lambda_{-j}=-\lambda_j,
\quad j=0,\ldots,m,\textrm{ and}
$$
$$
0=\lambda_0<\lambda_1<\ldots<\lambda_m<\pi.  
$$
Taking \begin{equation*}
\alpha=\min\{\alpha_j,\,\,j=0,\ldots,m\},\quad \text{and } \quad 
J=\{j \geq 0 \; : \; \alpha_j=\alpha\},  
\end{equation*}
 if $\alpha < \alpha_0/2$ then the condition \eqref{b(s)} is  satisfied.

Note that the condition on the coefficient $\alpha$ ensures
that $\sum_{h=1}^\infty |r(h)|=\infty $, thus the process has a long-memory.
However this condition is not enough to characterize the cyclical long
memory. 

\begin{enumerate}
\item When $\alpha<\alpha_0 /2$.  $\left \vert \displaystyle\sum_{j=1}^h r(j)
\right\vert = o\left(\displaystyle\sum_{j=1}^h r(j)^2 \right)$ as $h\to\infty$. 
Therefore the process  $(X_t^2)$ has  also a long
memory,  which is more persistent than $(X_t)$ (see Remark \ref{rem1}
for the exact expressions). 
This fact characterises cyclical long memory, and  the asymptotic behavior of many statistics (see
below for the empirical process) can be drastically different. 
  We focus on this case in this paper. 
\item When $\alpha > \alpha_0 /2$, the cyclical  behavior is less
  persistent than the regular long memory (singularity at frequency zero).  The
  presence of singularities outside zero do not modify  the convergence
  results obtained in the  regular case.  
\item When $\alpha = \alpha_0 /2$, both $(X_t)$ and $(X_t^2)$ will
  contribute to the limiting distribution, which will be a
  combinaison of the two
  previous cases.     
\end{enumerate}

We consider the empirical process associated with the
process $(X_n)_{n\ge1}$ defined by
$$
F_n(x)=\frac{1}{n}\sum_{j=1}^n\ind_{\{X_j\le x\}}.
$$

 \cite{ouldhaye:aphil:2003} proved the following results 
for the linear process $(X_n)$
defined in (\ref{g1}). Assume that $\BE
\xi_0^4<\infty$, the cumulative distribution function of $\xi_0$
is 5 times differentiable with continuous bounded and integrable
derivatives on $\mathbb{R}$. Denote
$$
d_n=n^{1-\alpha},\quad
 \textrm{and} \quad D=\frac{\sqrt{(2-2\alpha)(1-2\alpha)}}{4\Gamma(\alpha)\cos(\alpha\pi/2)}.
$$
If $\alpha<\alpha_0/2$, then, as $n$ tends
to infinity, we have
\begin{equation}\label{process empirique}
d_n^{-1}[nt]\bigl(F_{[nt]}(x)-F(x)\bigr)\Longrightarrow
\frac{F''(x)}{2}R(t) ,
\end{equation}
where $R$ is a linear combination of independent Rosenblatt
processes with the same parameter $\alpha$
\begin{equation}\label{equa9}
R(t)=R_{\alpha,\Lambda}(t)=D^{-1}\sum_{j\in J}c_j\Bigl(R_j^{(1)}(t)+R_j^{(2)}(t)\Bigr),
\end{equation}
where $\Lambda=\{\lambda_j,\quad j\in J\}$, and where 
\begin{itemize}
\item$c_0=h_0/2$, $c_j=h_j$ if $j\neq0$ and 
\begin{equation*}
h_j=g(e^{i\lambda_j})\prod_{k\neq
 j}\bigl(1-e^{i(\lambda_k-\lambda_j)}\bigr)^{(\alpha-1)/2},
\end{equation*} 
\item $R_j^{(i)}(t),\,i=1,2\textrm{ and }j\in J$
are Rosenblatt processes with parameter $1-\alpha$, independent
except for $j=0$, $R_0^{(1)}(t)=R_0^{(2)}(t)$.
\end{itemize}

The paper is organized as follows. In Section 2, we establish 
a limit theorem  for the kernel estimate. This extends one of
\cite{0862.60026}'s results, in particular we  show  the
contribution and the effect of the  singularities of the spectral
density  outside  the origin to the convergence rate and the limiting 
distribution. Then  we apply our limit theorem to construct
confidence bands for the density function.   

Similarly to \cite{0695.60043,MR1457496},  we provide  in Section 3,
the asymptotic behavior of the mean  integrated squared error,  and we
show that  the equivalence, one had in regular long memory' can be modified  when the  singularities of the
spectral density are not limited to the origin.

\section{Asymptotic distribution of the kernel estimator}\label{Sec:kernel}

Hereafter, we assume that the kernel $K$ 
is a continuous function with
compact support and $\int K(x) dx =1$. Concerning the bandwidth $m_n$, we assume that
 $m_n\to 0$ and $nm_n\to\infty$, as $n$ tends to infinity.\\
The equality
\begin{equation}\label{noyaux}
\tilde{f}_n(x)-\mathbb{E}\tilde{f}_n(x)=\frac{1}{m_n}\int_\mathbb{R}K\bigl(\frac{x-u}{m_n}\bigr)d\bigl(F_n(u)-F(u)\bigr)
\end{equation}
clearly shows the relationship between the estimate
$\tilde{f}_n(x)$ and the empirical process $F_n(x)$. The process
$\tilde{f}_n(x)$ is sometimes called the empirical density process.\\

For every integer  $n \ge1$,  we define the following statistics 
\begin{equation}\label{Y_{n,p}}
Y_{n,1}=\sum_{k=1}^n X_k,    \; \qquad
Y_{n,2}=\sum_{k=1}^n  \sum_{ s< r }   b_r b_{s}\xi_{k-s}
\xi_{k-r},
\end{equation}
and
\begin{equation}\label{S_{n,2}}
S_{n,2}(x)=n\bigl(F_n(x)-F(x))+F'(x)Y_{n,1}-\frac{1}{2}F''(x)Y_{n,2}.
\end{equation}

\begin{Rem}\label{rem1} For linear processes defined in \eqref{g1}, the following
equivalences as $n$ tends to infinity, have been proved by \cite{ouldhaye:aphil:2003}
 \begin{equation}\label{Y_{N,2}}
 \Var(Y_{n,2})\sim\frac 14 \Var\Bigl(\sum_{j=1}^n(X_j^2-\mathbb{E}(X_1^2))\Bigr)\sim
 Cn^{2-2\alpha}.
 \end{equation}
and 
\begin{equation}\label{var1}
\Var(Y_{n,1})=\Var\Bigl(\sum_{j=1}^nX_j\Bigr) \sim
Cn^{2-\alpha_0}. 
\end{equation}

Therefore (\ref{var1}) and (\ref{Y_{N,2}}) imply that the convergence rate
 obtained in Proposition \ref{estimation} is smaller than the convergence rate
of $\bar{X}_n$.
\end{Rem}

Let us define the class of Parzen kernels of order $s$.
\begin{Def}
  A kernel function $K$ is said to be a Parzen kernel of order $s\ge2$
  if it satisfies the following conditions
  \begin{enumerate}
  \item $ \int_{\mathbb{R}}K(u)du=1, $
  \item for every $1\le j\le s-1$, $\int_{\mathbb{R}}u^jK(u)du=0, $
  \item $\int_{\mathbb{R}}\vert u^s\vert\vert K(u)\vert du<\infty.$
  \end{enumerate}
\end{Def}
\cite{Bretagnolle} proved the existence of such kernels, for
which, an explicit construction  can be found in \cite{MR564251}.

\begin{Prop} \label{estimation}
Consider a process $(X_n)$
defined in (\ref{g1}). Assume that  $\alpha<\alpha_0/2$, $\BE
\xi_0^4<\infty$, the cumulative distribution function of $\xi_0$
is 5 times differentiable with continuous bounded and integrable
derivatives on $\mathbb{R}$.
 Let $K$ be a Parzen kernel of order $4$ having bounded total
variation. Assume that the bandwidth has the form
$$
m_n=n^{-\delta},\quad\textrm{where
}\frac{\alpha}{4}<\delta<\frac{\alpha}{2}.
$$
Then, as $n$ tends to infinity
\begin{equation}\label{HoSup}
n^\alpha\underset{x\in\mathbb{R}}{\sup}\vert
\tilde{f}_n(x)-f(x)\vert
\overset{d}{\longrightarrow}\underset{x\in\mathbb{R}}{\sup}
\Bigl\vert\frac{f''(x)}{2}\Bigr\vert\vert R_{\alpha,\Lambda}\vert.
\end{equation}

where $R_{\alpha,\Lambda}=R_{\alpha,\Lambda}(1)$.
Moreover,
\begin{equation}\label{vectorielle}
n^{\alpha}(\tilde{f}_n(x)-f(x))\overset{C_b(\mathbb{R})}{\Longrightarrow}
-\frac{f''(x)}{2}R_{\alpha,\Lambda},
\end{equation}
where $\overset{C_b(\mathbb{R})}{\Longrightarrow}$ denotes the
convergence in $C_b(\mathbb{R})$, the space of continuous bounded
functions.
\end{Prop}

Proof:\\
The difference between $\tilde{f}_n$ and $f$ can be expressed as
\begin{align*}
\tilde{f}_n(x)-f(x)& =\tilde{f}_n(x)-\mathbb{E}\tilde{f}_n(x)+\mathbb{E}\tilde{f}_n(x)-f(x)
\\
&=\frac{1}{m_n}\int K(u)d\bigl(F_n(x-m_nu)-F(x-m_nu)\bigr)
+\int\bigl(f(x-m_nu)-f(x)\bigr)K(u)du.
\end{align*}

We first replace $F_n-F$ by its expression in (\ref{S_{n,2}}).
Then we
apply the integration by parts formula on the first integral. For the
second, we apply the Taylor-Lagrange formula. There exists
a real number $u^*$ such that $\vert u^*-x\vert<\vert m_nu\vert$ and
\begin{align*}
\tilde{f}_n(x)-f(x)=& \frac{-1}{nm_n}\int
S_{n,2}(x-m_nu)dK(u)+\frac{Y_{n,1}}{n}\int
f'(x-m_nu)K(u)du\\
&-\frac{Y_{n,2}}{n}f''(x)\int K(u)du+\frac{Y_{n,2}}{n}m_n\int
f^{(3)}(u^*)uK(u)+\\
&+\int
\bigl(-m_nuf'(x)+\frac{m_n^2u^2}{2}f''(x)-\frac{m_n^3u^3}{6}f^{(3)}(x)+
\frac{m_n^4u^4}{24}f^{(4)}(u^*)\bigr)K(u)du\\
=:&a_n(x)+b_n(x)+c_n(x)+d_n(x)+e_{n}(x).
\end{align*}

 Now, a proof similar to that of Theorem 2.2 in \cite{0862.60026} 
allows us to write for $2\delta<\alpha$
\begin{equation}\label{ho}
n^{\alpha+\delta-1}\underset{x\in\mathbb{R}}{\sup}\vert
S_{n,2}(x)\vert\overset{a.s.}\longrightarrow 0, \qquad \textrm{as} \;
n\to \infty.
\end{equation}
And, thus we have
\begin{equation}\label{probability}
n^\alpha\underset{x\in\mathbb{R}}{\sup}\vert a_n(x)\vert
\overset{P}\longrightarrow0, \qquad \textrm{as} \;
n\to \infty
\end{equation}
where $\overset{P}\longrightarrow$ denotes the convergence in probability.\\
For the sequences $b_n(x)$, $d_n(x)$, $e_n(x)$, we get the
same convergence in probability as in (\ref{probability}) by bounding the
variances. To obtain the bounds, we start from the variances of $Y_{n,1}$
and $Y_{n,2}$ defined in (\ref{var1}) and (\ref{Y_{N,2}}), and we
use the fact that $K$ is a Parzen kernel and that $f$ is 4 times
differentiable and bounded derivatives. We get, as $n$ tends to
infinity,
\begin{eqnarray*}
\Var(n^\alpha\underset{x\in\mathbb{R}}{\sup}\vert
b_n(x)\vert)&\le&
n^{2\alpha-2}\Var\Bigl(Y_{n,1}\underset{x\in\mathbb{R}}{\sup}\vert
f'(x)\vert\int\vert K(u)\vert du\Bigr)\\
&=&Cn^{2\alpha-2}n^{2-\alpha_0}=Cn^{2\alpha-\alpha_0}\longrightarrow 0,
\end{eqnarray*}
\begin{eqnarray*}
\Var(n^\alpha\underset{x\in\mathbb{R}}{\sup}\vert
d_n(x)\vert)&\le&
n^{2\alpha-2}\Var\Bigl(Y_{n,2}m_n\underset{x\in\mathbb{R}}{\sup}\vert
f^{(3)}(x)\vert\int\vert u K(u)\vert du\Bigr)\\
&=&Cn^{2\alpha-2}n^{2-2\alpha}n^{-\delta}\longrightarrow 0,
\end{eqnarray*}
\begin{equation*}
n^\alpha\underset{x\in\mathbb{R}}{\sup}\vert
e_n(x)\vert\le\underset{x\in\mathbb{R}}{\sup}\vert
f^{(4)}(x)\vert\frac{n^{\alpha-4\delta}}{24}\int u^4\vert
K(u)\vert du=O(n^{\alpha-4\delta})\longrightarrow 0,
\end{equation*}
These four convergences in probability imply that both sequences
$$
n^\alpha\underset{x\in\mathbb{R}}{\sup}\vert\tilde{f}_n(x)-f(x)\vert\quad\textrm
{and}\quad n^\alpha\underset{x\in\mathbb{R}}{\sup}\vert f''(x)\vert
\big\vert\frac{Y_{n,2}}{n}\big\vert=n^\alpha\underset{x\in\mathbb{R}}{\sup}\vert
c_n(x)\vert
$$
have the same limit as $n$ tends to infinity. According to Lemma~2.1
in \cite{ouldhaye:aphil:2003}, this common limit is equal to
$$ \underset{x\in\mathbb{R}}{\sup}
\Bigl\vert\frac{f''(x)}{2}\Bigr\vert\vert R_{\alpha,\Lambda}\vert.
$$
 Hence (\ref{HoSup}) is proved. According to (\ref{Y_{N,2}}), we notice that the rate
 $n^{-\alpha}$ given in (\ref{HoSup}) is the convergence rate of
 $n^{-1}\sum_{j=1}^n(X_j^2-\mathbb{E}(X_1^2))$.\\
Similarly, as $n$ tends to infinity, the finite-dimensional
distributions of
$$
n^{\alpha}(\tilde{f}_n(x)-f(x))\quad\textrm {and}\quad-n^\alpha
f''(x)\frac{Y_{n,2}}{n}=n^\alpha c_n(x)
$$
converge simultaneously to the finite-dimensional distributions of
$ -(f''(x)/2)R_{\alpha,\Lambda}. $ This concludes the proof of
(\ref{vectorielle}) because (\ref{HoSup})
implies the tightness of $n^{\alpha}(\tilde{f}_n(x)-f(x))$.

\begin{Rem} We clearly see that the choice of the class
 of Parzen kernels allows the bias
$\mathbb{E}\tilde{f}_n(x)-f(x)$ to become negligible. If $K$ is not
a Parzen kernel, the contribution of the bias $e_n(x)$ is
not negligible with respect to $b_n(x)$. Therefore, (\ref{HoSup})
is false for a standard kernel unless we replace
$\tilde{f}_n(x)-f(x)$ by $\tilde{f}_n(x)-\mathbb{E}\tilde{f}_n(x)$
in (\ref{HoSup}).
\end{Rem}

\begin{Rem}
The result (\ref{HoSup}) in Proposition
\ref{estimation} can be applied to obtain a goodness of fit test on
the marginal density.
\end{Rem}

\begin{Rem} The result (\ref{HoSup}) in Proposition
\ref{estimation} provides confidence bands for $f$ which
depend on the derivative $f''$. In general, $f''$ is not
available, and thus the confidence band cannot be calculated. Then
 $f''$ can be replaced by its kernel estimate given by
$$
\tilde{f}_n''(x)=\frac{1}{nm_n^3}\sum_{j=1}^nK''\bigl(\frac{x-X_j}{m_n}\bigr).
$$
(note that it is necessary to assume that the kernel function $K$
is twice differentiable.)
\end{Rem}

\begin{Prop}\label{confidence}
Under the same hypotheses as in Proposition \ref{estimation} and if
the kernel function $K$ is twice differentiable and its
derivative $K''$ is continuous, then for each interval
$[a,b]$ on which $f''$ is positive, we have
\begin{equation}\label{confiance}
2n^\alpha\underset{x\in[a,b]}{\sup}\Big\vert\frac{\tilde{f}_n(x)-f(x)}{\tilde{f}_n''(x)}
\Big\vert\overset{d}{\longrightarrow}\vert R_{\alpha,\Lambda}\vert.
\end{equation}
\end{Prop}
In other words, as $n$ tends to infinity, for every $t>0$, we have
\begin{equation}
P\big\{\tilde{f}_n(x)-\frac{t\tilde{f}_n''(x)}{2n^\alpha}\le
f(x)\le\tilde{f}_n(x)+\frac{t\tilde{f}_n''(x)}{2n^\alpha},\,\,a\le
x\le b\big\}\to P\big\{\vert R_{\alpha,\Lambda}\vert<t\big\}.\label{eq:15}
\end{equation}

In  Proposition \ref{prop-quantil}, we give a  consistent estimate of
 the quantiles of process  $R_{\alpha,\Lambda} $.   Using
(\ref{eq:15}), this allows us
to  obtain  asymptotic confidence band for the density
$f(x)$ which is valid for every $x\in[a,b]$.\\

\noindent \textbf{Proof : } \\ 
Let $\phi$ be the function defined on $C_b(\mathbb{R})$ by
$$\phi(g)=\underset{x\in[a,b]}{\sup}\Big\vert\frac{g(x)}{f''(x)}\Big\vert $$
Since $\phi$ is continuous,   (\ref{vectorielle})  ensures the
following convergence : 
\begin{equation}\label{continue}
2n^\alpha\underset{x\in[a,b]}{\sup}\Big\vert
\frac{\tilde{f}_n(x)-f(x)}{f''(x)}\Big\vert
\overset{d}{\longrightarrow}\vert R_{\alpha,\Lambda}\vert, \qquad \textrm{as}\; n\to\infty.
\end{equation}
Now, we prove that the difference
$$
Y_n(x):=n^\alpha\Bigl(\frac{\tilde{f}_n(x)-f(x)}{f''(x)}
-\frac{\tilde{f}_n(x)-f(x)}{\tilde{f}_n''(x)}\Bigr)
$$
satisfies
$$
\underset{x\in[a,b]}{\sup}\vert
Y_n(x)\vert\overset{P}{\longrightarrow}0, \qquad \textrm{as}\; n\to\infty.
$$
This convergence is obtained as follows. We rewrite $Y_n(x)$ as
$$
\vert
Y_n(x)\vert=n^\alpha\Big\vert\frac{\tilde{f}_n(x)-f(x)}{f''(x)}\Big\vert
\Big\vert\frac{\tilde{f}_n''(x)-f''(x)}{\tilde{f}_n''(x)}\Big\vert.
$$
and by (\ref{continue}), it is enough to prove that
\begin{equation}\label{continue1}
\underset{x\in\mathbb{R}}{\sup}\Big\vert\frac{\tilde{f}_n''(x)-f''(x)}{\tilde{f}_n''(x)}\Big\vert
\overset{P}{\longrightarrow}0, \qquad \textrm{as}\; n\to\infty.
\end{equation}
The difference between
$\tilde{f}_n''$ and $f''$ can be written as
\begin{align*}
\tilde{f}_n''(x)-f''(x) = &  \frac{-1}{nm_n^3}\int
S_{n,2}(x-m_nu)dK''(u)+\frac{Y_{n,1}}{n}
\int f^{(3)}(x-m_nu)K(u)du - \\
&-\frac{Y_{n,2}}{n}\int f^{(4)}(x-m_nu)K(u)du
+\int\bigl(f''(x-hu)-f''(x)\bigr)K(u)du.
\end{align*}
 by replacing $f$ with $f''$ and $\tilde{f}_n$ with $\tilde{f_n''}$
and following the same lines as the proof of Proposition \ref{estimation}. Then, we get
$$
\underset{x\in\mathbb{R}}{\sup}\vert\tilde{f}_n''(x)-f''(x)\vert
=O\bigl(n^{-(2\delta\wedge(1-3\delta))}\bigr).
$$
Since $0<\delta<1/4$, we have
$$
\underset{x\in\mathbb{R}}{\sup}\vert\tilde{f}_n''(x)-f''(x)\vert\overset{P}{\longrightarrow}0,
$$
moreover, the derivative $f''$ satisfies
$$
\underset{x\in[a,b]}{\inf}\vert f''(x)\vert>0.
$$
Thus, we get (\ref{continue1}). This concludes the proof.\\

\begin{Prop}\label{prop-quantil}
 Fix $\beta\in(0,1)$. Let $c(\alpha,\Lambda,\beta)$ be the quantile
 of order $\beta$ of the process $R_{\alpha,\Lambda}$ defined in
 \eqref{equa9}. 
 If $(\alpha_n,\Lambda_n)$ are consistent (in probability) estimators of
 $(\alpha,\Lambda)$.
then
\begin{equation}
 c(\alpha_n,\Lambda_n,\beta)\overset{P}{\to}c(\alpha,\Lambda,\beta)\label{eq:1}
\end{equation}

\end{Prop}
\begin{Rem}
In the  references given in the introduction, the parametric and semi
parametric methods provide estimators of  $(\alpha,\Lambda)$ which satisfy the
condition required in Proposition \ref{prop-quantil}. 
\end{Rem} 
{Proof : }
We want to show \eqref{eq:1} which will be obtained if we show that the application
$(\gamma,\theta) \mapsto c(\gamma,\theta,\beta)$ is continuous, as $(\alpha_n,\Lambda_n)\overset{P}{\to}(\alpha,\Lambda)$. To prove this continuity we prove that the mappings $g,h$ below are continuous,
$$
((0,1)\times [0,\pi]^{|J|}, \;
\vert.\vert)\overset{g}{\to}(C_b(\mathbb{R}),\|.\|)\overset{h}{\to}((0,1),\vert.\vert),$$
where 
$\|.\|$ is the uniform metric, and in the following decomposition
$F_{\gamma,\theta}$ is the distribution function of
$R_{\gamma,\theta}$.
$$
(\gamma,\theta) \mapsto [g(\gamma,\theta)=F_{\gamma,\theta}]\mapsto
[h(F_{\gamma,\theta})=c(\gamma,\theta,\beta)].$$
Continuity of $g$ can be proved as follows. Consider a deterministic
sequence $(\gamma_n,\theta_n)$ such that
$(\gamma_n,\theta_n)\to(\gamma,\theta)$ as $n\to\infty$. Then to prove that
$F_{\gamma_n,\theta_n}\to F_{\gamma,\theta}$ uniformly it will be
enough to show that $R_{\gamma_n,\theta_n}\Longrightarrow
R_{\gamma,\theta}$. To obtain the latter weak convergence it will
suffice to show that every sequence of Rosenblatt variables
$(R_{\gamma_n})$ with parameter $\gamma_n$ converges weakly to a
Rosenblatt variable $R_{\gamma}$ with parameter $\gamma$, as
$R_{\gamma_n,\theta_n}$ is a linear combination of independent
Rosenblatt variables $R_{\gamma_n}$ with the coefficients $c_j/D$ that
are continuous functions of $\gamma_n,\theta_n$. We have from
\cite{major} $$
R_{\gamma_n}=\int\int_{\mathbb{R}^2}\frac{e^{i(x+y)}-1}{i(x+y)}W_n(dx,dy)
$$
where 
$$W_n(dx,dy)= \vert x\vert^{(\gamma_n-1)/2}\vert y\vert^{(\gamma_n-1)/2}W(dx,dy)
$$ 
with $W(dx,dy)$ being the standard Gaussian random measure, and since

$$\vert x\vert^{(\gamma_n-1)/2}\vert y\vert^{(\gamma_n-1)/2}\to\vert x\vert^{(\gamma-1)/2}\vert y\vert^{(\gamma-1)/2}
$$ then we have the required convergence. \\
Now to prove the continuity of $h$ it is enough to note that the
quantile function is continuous (with respect to the uniform metric)
over the class of monotonic continuous distribution functions, i.e. if
$\|F_n-F\|\to0$ then $h(F_n,\beta))\to h(F,\beta)$. Of course here we
do have 
 $\|F_{\gamma_n,\theta_n}-F_{\gamma,\theta}\|\to0$, as we
just established that $R_{\gamma_n,\theta_m}\Longrightarrow
R_{\gamma,\theta}$.\hfill $\square$

\section{Asymptotic mean integrated squared error (MISE)}

The mean integrated squared error (MISE) of the estimate
$\tilde{f}_n$ is defined by 
$$
\int_\mathbb{R}\mathbb{E}\bigl(\tilde{f}_n(x)-f(x)\bigr)^2dx.
$$
For a wide class of linear processes including the processes
with short and \textit{regular} long memories,
\cite{0695.60043} and \cite{MR1457496} studied the asymptotic
behavior of the MISE. In particular, they established the following equivalence, when $n$ tends to infinity,
\begin{equation}\label{Hall1}
\int_\mathbb{R}\mathbb{E}\bigl(\tilde{f}_n(x)-f(x)\bigr)^2dx\sim
 \int_\mathbb{R}\mathbb{E}_0\bigl(\tilde{f}_n(x)-f(x)\bigr)^2dx
+\Var(\overline{X}_n)\int_{\mathbb{R}}f'(x)^2dx
\end{equation}
where $\BE_0$ denotes the expectation with respect to the distribution of
$n$ independent random variables distributed from the density $f$.
In particular, the equivalence (\ref{Hall1}) shows that the
convergence rate of the MISE cannot be faster than
the convergence of $\Var(\overline{X}_n)$. In
other words, the convergence rate of the kernel density estimates is
bounded from above by the convergence rate of the empirical mean.
This is the optimal rate.

Hereafter, we assume that the distribution of the innovation $(\xi_k)$ satisfies 
\begin{itemize}
\item[\textbf{[Z]}] There exist $\delta > 0$ and $C<\infty$  such that
  the characteristic
 function of $\xi_0$ satisfies 
\begin{equation}
|E e^{iu\xi_0}| \leq C (1+|u|)^{-\delta}\label{eq:z2}
\end{equation}

 \end{itemize}

\begin{Theo}\label{sec:asympt-mean-integr}
 Let $(X_n)$ be a linear process 
defined in (\ref{lineaires}) and (\ref{g1}) such that  the
distribution of $\xi_0$ satisfies $[Z]$ and 
 $\BE
\xi_0^4<\infty$.
Assume that $\alpha<\frac 13 \wedge\frac{\alpha_0}{2}$ and the kernel $K$ is a
bounded symmetric density function. 
Then the MISE satisfies, as $n$ tends to infinity, 
\begin{align}
 MISE(\tilde f_n) \sim  &
 \int_\mathbb{R}\mathbb{E}_0\bigl(\tilde{f}_n(x)-f(x)\bigr)^2dx
+\frac 14 \Var \Bigl(\frac 1n \sum_{j=1}^n(X_j^2-\mathbb{E}(X_1^2))\Bigr) \int_{\mathbb{R}}f''(x)^2dx
\label{eq:2}
\end{align}
where $\BE_0$ denotes the expectation with respect to the distribution of
$n$ independent random variables distributed from the density $f$.
\end{Theo}

\begin{Rem}
The variance $ \Var \Bigl(\frac 1n \sum_{j=1}^n(X_j^2-\mathbb{E}(X_1^2))\Bigr) $ is also equivalent to 
$4\Var(\frac 1n {Y}_{n,2}) $ (see \cite{ouldhaye:aphil:2003}). 
Equation \eqref{eq:2} shows that this term is
 a ceiling rate of MISE independently of the choice of the kernel and
 bandwidth. 
\end{Rem}

 \noindent{\bf Proof. }~ \\
\textit{Notation :} for an  arbitrary function $g$,  we denote by $\hat g $ its  Fourier
transform. 

The proof consists in adapting  the proof of \cite{MR1457496} to the cyclical
case.
Using \cite{0695.60043} decomposition of the MISE, we have 
\begin{align}
 \mathrm{MISE}(\tilde f_n) = &
 \int_\mathbb{R}\mathbb{E}_0\bigl(\tilde{f}_n(x)-f(x)\bigr)^2\d x +  \nonumber\\ 
&+\frac1{n\pi} \sum_{j=1}^{n-1} (1-j/n) \int |\hat K (m_nt)|^2 \left\{ \Re(
 \mathbb{E} (e^{i t (X_1-X_{j+1})} ) - |\hat f (t)|^2\right\} \d t \label{eq:5}\\ 
& :=
\mathrm{MISE}_0 + W_n \nonumber
\end{align}
Let $f_j$ be the joint density of $(X_1,X_{j+1})$. 
We extend the expansion of $f_j$, obtained by 
\cite{MR1409327} for the first order, to the second order  as follows: there exists a function
$\ell_j \; : \mathbb{R}^2 \mapsto \mathbb{R} $ such
that  
\begin{equation}
f_j(x,y) = f(x) f(y) +r(j) f'(x) f'(y) + \frac 12 r(j)^2 f''(x)
f'' (y) +\ell_j(x,y) \quad \forall (x,y)\in \mathbb{R}^2\label{eq:3}
\end{equation}
where $r$ is given in \eqref{r(n)}. 

We have 
\begin{align}
 \label{eq:4}
 \mathbb{E} (e^{i t (X_1-X_{j+1})} ) & = \int e^{it(x-y) } f(x) f(y)
 \d x \d y +r(j) \int e^{it(x-y) } f'(x) f'(y)\d x \d y +\nonumber \\ & +\frac 12 r(j)^2\int
 e^{it(x-y) } f''(x)
f'' (y) \d x \d y + \int e^{it(x-y) } \ell_j(x,y)\d x \d y \nonumber \\
 &= |\hat f (t)|^2 +r(j) |\widehat{f'}(t)|^2 +\frac 12 r(j)^2
 |\widehat{f''}(t)|^2 +\widehat{ \ell_j}(t,-t). 
\end{align}
Similarly to \cite{MR1457496}, $W_n$ in \eqref{eq:5} can be written as 
\begin{multline*}
 W_n = \frac{2}{n} \sum_{j=1}^{n-1} (1-j/n) r(j) \int |K_{m_n}\star
 f'|^2(t) \d t + \frac{1}{n} \sum_{j=1}^{n-1} (1-j/n) r(j)^2 \int
 |K_{m_n}\star f''|^2(t) \d t +\\+\frac{1}{n\pi} \sum_{j=1}^{n-1}
 (1-j/n) \int |\hat K(m_n t)|^2 \Re \widehat{\ell_j}(t,-t) \d t
\end{multline*}
where $K_{m_n}(x) = m_n^{-1}K(x m_n^{-1})$, and  where  $f\star g$  is
the convolution of $f$ and $g$. Moreover  we have,  for $k = 1,2$, 
$$\int |K_{m_n}\star f^{(k)}|^2(t) \d t = \int f^{(k)}(t)^2 \d t +
o(1), \quad n\to\infty.$$ 

We obtain 
\begin{multline}\label{eq:10}
 W_n = \frac2 {n} \sum_{j=1}^{n-1} (1-\frac jn) r(j) \left( \int f'(t)^2 \d
 t +o(1) \right) + \frac1 {n} \sum_{j=1}^{n-1} (1-\frac jn)
r(j)^2\left( \int f''(t)^2 \d t +o(1) \right) + \\ + \frac{1}{n\pi}
 \sum_{j=1}^{n-1} (1-\frac jn) \int |\hat K(m_n t)|^2 \Re \widehat{\ell_j}(t,-t)
 \d t.
\end{multline}
According to \cite{MR1027992}, we have 
\begin{align}
 \Var \Bigl(\frac 1n
 \sum_{j=1}^n(X_j^2-\mathbb{E}(X_1^2))\Bigr) & = \frac{2}{n^2}
 \sum_{1\leq i,j\leq n } r^2(i-j) +O(n^{-1}), \nonumber \\ 
 &=\frac2 {n^2} (nr(0) + 2 \sum_{j=1}^{n-1} (n-j) r(j)^2 )
+O(n^{-1}) \nonumber \\ & = \frac 4 {n} \sum_{j=1}^{n-1} (1-j/n) 
r(j)^2 +O(n^{-1}) :=\gamma(n) \label{eq:11}
\end{align}
Moreover, using the form of $r$ given in \eqref{r(n)} and the fact
that $\alpha < 1/3$, we get 
\begin{align}
\gamma(n) &= \frac 4 {n} \sum_{j=1}^{n-1} (1-j/n) j^{-2\alpha}
(\sum_{k\in J}a_k \bigl(\cos j \lambda_k+o(1)\bigr))^2 
 +O(n^{-1}) \nonumber \\
&= \frac 2 {n} \sum_{j=1}^{n-1} (1-j/n)
j^{-2\alpha} \sum_{k\in J} a_k^2
+O(n^{-1}) \nonumber \\
&= \frac 2 {n} n^{1-2\alpha} \left( \frac{1}{1-2\alpha} -\frac{1}{2-2\alpha}\right) \sum_{k\in J} a_k^2
+O(n^{-1}) \nonumber \\ 
&= n^{-2\alpha} \frac{1}{(1-2\alpha)(1-\alpha)} \sum_{k\in J} a_k^2
+O(n^{-1}) \sim C n^{-2\alpha} \label{eq:7}
\end{align}
As  $\alpha < \frac 13 \wedge \frac{\alpha_0}{2}$ and using \eqref{var1}, we get 
\begin{equation}
\frac2{n} \sum_{j=1}^{n-1} (1-j/n) r(j) =\frac 1{n^2}\Var(Y_{n,1})
-r(0) n^{-1}=  O(n^{-\alpha_0}) + O(n^{-1})=
o(n^{-2\alpha}).\label{eq:6}
\end{equation}

From \eqref{eq:10}, \eqref{eq:11}, \eqref{eq:7} and \eqref{eq:6}, we get 
\begin{align*}
 W_n =\frac 14 \Var \Bigl(\frac 1n
\sum_{j=1}^n(X_j^2-\mathbb{E}(X_1^2))\Bigr) &\int f''(t)^2 \d t
+o(n^{-2\alpha}) + \\ &+ \frac{1}{n\pi}
 \sum_{j=1}^{n-1} (1-j/n) \int |\hat K(m_n t)|^2 \Re \widehat{\ell_j}(t,-t)
 \d t.
 \end{align*}
 Since $r(j)^2$ behaves asymptotically as $j^{-2\alpha}$, and 
\begin{equation}
\frac{1}{n\pi}
 \sum_{j=1}^{n-1} (1-j/n) \int |\hat K(m_n t)|^2 \Re \widehat{\ell_j}(t,-t)
 \d t \leq \frac{1}{n\pi}
 \sum_{j=1}^{n-1} (1-j/n) \int | \widehat{\ell_j}(t,-t)|
 \d t \end{equation}
the proof is completed using the following lemma proven below. 
\begin{Lem}\label{sec:asympt-mean-integr-1}
 Under the same assumption of Theorem \ref{sec:asympt-mean-integr},
\begin{equation}
\int |\widehat{\ell_j}(t,-t)| 
 \d t = O(j^{-2\alpha - \epsilon} ). \label{eq:9}
\end{equation}
for $\epsilon$ an arbitrary positive number smaller than
$  \dfrac{1-3\alpha}{10}$.  

\end{Lem}

\null \hfill $\square$

\noindent \textbf{Proof of Lemma \ref{sec:asympt-mean-integr-1}}
By definition of $\ell_j$ in \eqref{eq:3}, we have 
$$
\widehat{\ell_j}(x,y) = \widehat{f_j}(x,y) - \widehat f(x) \widehat f(y)(
1-xy r(j) +\frac 12 x^2y^2 r(j)^2) $$

We split the integral 
\begin{equation}\label{eq:8}
 \int_{\mathbb{R}} |\widehat{\ell_j}(t,-t)| 
 \d t = \int_{|t|> j^\epsilon}
 |\widehat{\ell_j}(t,-t)| \d t + \int_{|t|<j^\epsilon} |\widehat{\ell_j}(t,-t)| \d t
\end{equation}
where  $\epsilon$ is an arbitrary positive number smaller than $  \dfrac{1-3\alpha}{10}$.  

Under assumption \eqref{eq:z2}, 
\cite{MR1409327} proved for the regular long memory that for arbitrary
$k$ 
$$|\widehat{f_j}(x_1,x_2) | \leq c(k) (1+|x|)^{-k}$$  for all
$x=(x_1,x_2)\in\mathbb{R}^2$ and 
$$|\widehat{f}(x) | \leq c(k) (1+|x|)^{-k}$$ for all $x\in\mathbb{R}$. 

Their proof can be adapted to the cyclical case i.e. when the coefficients
$(b_j)_{j\in\mathbb{N}}$ satisfies \eqref{b(s)}. Using their
notation, it suffices to construct a finite set $J_1$ such that for
all $j\in J_1$ : $|b_{-j} |> 2| b_{t-j}| +c_1$ where $c_1$ does not
depend on $t$. 
Since $(|b_j|)_{j\in\mathbb{Z}}$ is not summable, there exists a subsequence
$({j_u})_{u\in \mathbb{Z}}$ such that $b_{-j_u} \not = 0$. We can take
$J_1 $ a subset of $\{ j_u : u\in\mathbb{Z} \}$ with $[\delta
|J_1|] =k+3$. 
Indeed, for $j\in J_1$, we have 
$|b_{-j}| > C(J_1) |j^{-(\alpha+1)/2}| $, and for $t$ large enough
there  exists $\tilde c_1$ 
 $$ |j|^{-(1+\alpha)/2} > 2 /C(J_1)
 |t-j|^{-(1+\alpha)/2} +\tilde c_1.$$ 
Therefore, there exists $c_1$ such that for all $j\in J_1$, 
$$ |b_{-j} | > 2| b_{t-j}| + c_1.$$

For all  $k'$, the first integral in (\ref{eq:8}) satisfies 
$$ \int_{|t|> j^\epsilon}
 |\widehat{\ell_j}(t,-t)| \d t \leq  j^{-\epsilon k'} \int_{|t|> j^\epsilon} |t|^{k'}
|\widehat{\ell_j}(t,-t)| \d t = O(j^{-\epsilon k'}).$$ 
Therefore we can take any arbitrary $k'$ such that $k'> (2\alpha+\epsilon)/\epsilon$. 

For the second integral in \eqref{eq:8}, it is enough to show that 
\begin{equation}
\sup_{|u|< j^\epsilon} |\widehat{\ell_j}(u) | =
O(j^{-2\alpha-2\epsilon }).\label{eq:12}
\end{equation}

The proof is quite similar to that of equation $(2.20)$ in Giraitis \textit{et al}
\cite{MR1409327} adding the terms of order two in the
expansion. 

We write the difference $\widehat{f_j}(x,y) - \widehat f(x) \widehat f(y)$
from products of the characteristic function $\phi$ of $\xi_1$. 

\begin{align*}
\widehat{f_j}(x,y) - \widehat f(x) \widehat f(y) & = \prod_{I_1}
\prod_{I_1} \prod_{I_1} \phi(xb_{-i} +y b_{t-i}) - \prod_{I_1}
\prod_{I_1} \prod_{I_1} \phi(xb_{-i} )\phi(y b_{t-i}) := a_1a_2a_3
-a'_1a'_2a'_3\\
&= (a'_1-a_1) a_2a_3 + (a'_2-a_2) a'_1a_3+ (a'_3-a_3)a'_1a'_2 \end{align*}
where $I_1= \{ |i| < j^{2\epsilon}\}$, $I_3= \{ |t-i| < j^{2\epsilon}\}$
and $I_3= \mathbb{Z}-(I_1\cup I_2)$. 
We will deduce \eqref{eq:12} from $|a_i|<1$, $|a'_i|<1$ and the following
facts, for all $u<t^\epsilon$ 
\begin{eqnarray}
 \label{eq:13}
 a_i-a'_i &=& O(j^{-2\alpha-2\epsilon}), \quad i=1,2\\
a_3-a'_3 &= & a'_3(-xy r(j) +\frac 12 x^2y^2 r(j)^2) + O(j^{-2\alpha-2\epsilon}).\label{eq:14}
\end{eqnarray}

Similarly to \cite{MR1409327}, we prove \eqref{eq:13} with $i=1$ (or
similarly for $i=2$) as follows 
\begin{align*}
 | a_1-a'_1| &\leq \sum_{|i|\leq j^{2\epsilon}} |\phi(xb_{-i} +y
 b_{j-i}) - \phi(xb_{-i} )\phi(y b_{j-i}) |\\ 
&\leq \sum_{|i|\leq j^{2\epsilon}} |xb_{j-i} |
\end{align*}
As $|i|\leq j^{2\epsilon}$ and $x\leq j^\epsilon$, we have 
$$|xb_{j-i}|\leq C j^\epsilon j^{-(1+\alpha)/2} =
j^{-2\alpha-4\epsilon} O(1) $$ 
since $\epsilon<   \dfrac{1-3\alpha}{10}$. Therefore $|a_1 - a'_1| =
j^{-2\alpha-2\epsilon } O(1) .$ 

To prove \eqref{eq:14}, we follow the same calculations as \cite{MR1409327} page
325. Since $|xb_{-i} | + |yb_{j-i}| =o(1)$, we write $a_3-a'_3$ of the
form $$ a_3-a'_3 = a'_3 (e^{Q_j(x,y)}-1) = a'_3 (Q_j(x,y) + \frac 12
Q_j(x,y)^2 +o(Q_j(x,y)^2) ) $$ 
where 
\begin{align*}
  Q_j(x,y) &=\sum_{i\in I_3} \Psi(xb_{-i},yb_{j-i}) = -xy \sum_{i\in
    I_3} b_{-i} b_{j-i} +O(\sum_{i\in I_3} (x b_{-i})^2 |y b_{j-i}|
  +|x b_{-i}| |yb_{j-i}|^2 ) \\ &:= -xy \sum_{i\in I_3} b_{-i} b_{j-i}
  +R_n
\end{align*}
and $$ \Psi(x,y) = \log(\phi(x+y) ) -\log(\phi(x) )-\log(\phi(y) )$$ and we show that 
$$Q_j(x,y) =-xyr(j) + O(\sum_{I_1\cup I_2} |x||y||b_{-i}||b_{j-i}|
+\sum_i x^2 |y| |b_{-i}|^2|b_{j-i}| ) = $$
\begin{multline*}
 Q_j(x,y)^2 = x^2 y^2 r(j)^2+ x^2y^2 (\sum_{i\in I_1•\cup I_2} b_{-i}
 b_{j-i} )^2 - 2x^2y^2 -xy \sum_{i\in \mathbb{Z}} b_{-i}
 b_{j-i}\sum_{i\in I_1\cup I2} b_{-i} b_{j-i} \\ + R_n^2 -2R_n xy
 \sum_{i\in I_3} b_{-i} b_{j-i}
\end{multline*}
For $|x|< j^\epsilon$ et $|y|< j^\epsilon$ we have 
$$\sum_{I_1\cup I_2} |x||y||b_{-i}||b_{j-i}| =
j^{2\epsilon-(1+\alpha)/2 } O(1) = j^{-2\alpha-2\epsilon} O(1) $$ 
since $\epsilon < (1-3\alpha )/8 $ 
and 
$$\sum_i x^2 |y| |b_{-i}|^2|b_{j-i}| = j^{-\alpha/2 -1/2 +3\epsilon}
O(1) = j^{-2\alpha-2\epsilon} O(1) $$ 
since $\epsilon < (1-3\alpha )/10 $. 
These asymptotic behaviors ensure that for $|x|< j^\epsilon$ et $|y|< j^\epsilon$ we have 
$$ a_3-a'_3 = a'_3 ( xy r(j) +\frac 12 x^2 y^2 r(j)^2
+O(j^{-2\alpha-2\epsilon}) ). $$

\begin{center}
  \textbf{Acknowledgement}
\end{center}

The authors would like to thank the anonymous referee for their
helpful comments and suggestions, that  improved the presentation
of the paper.

\end{document}